\documentclass[10pt]{article}
\usepackage{bbm}
\usepackage{amsfonts}
\usepackage{amsmath}
\allowdisplaybreaks
\usepackage{amssymb}
\usepackage{fancyhdr}
\usepackage{setspace}
\usepackage{latexsym}
\usepackage{mathrsfs}
\usepackage{amsthm}

\usepackage{cite,enumitem,graphicx}
\usepackage[colorlinks=true,urlcolor=blue,
citecolor=red,linkcolor=blue,linktocpage,pdfpagelabels,
bookmarksnumbered,bookmarksopen]{hyperref}
\usepackage[english]{babel}

\numberwithin{equation}{section}
\newtheorem{theorem}{Theorem}[section]
\newtheorem{proposition}[theorem]{Proposition}
\newtheorem{lemma}[theorem]{Lemma}
\newtheorem{corollary}[theorem]{Corollary}
\newtheorem{definition}[theorem]{Definition}
\newtheorem{remark}[theorem]{Remark}
\newtheorem{example}[theorem]{Example}

\newcommand{\bt}{\begin{theorem}}
\newcommand{\et}{\end{theorem}}
\newcommand{\bl}{\begin{lemma}}
\newcommand{\el}{\end{lemma}}
\newcommand{\bd}{\begin{definition}}
\newcommand{\ed}{\end{definition}}
\newcommand{\bc}{\begin{corollary}}
\newcommand{\ec}{\end{corollary}}
\newcommand{\bp}{\begin{proof}}
\newcommand{\ep}{\end{proof}}
\newcommand{\bx}{\begin{example}}
\newcommand{\ex}{\end{example}}
\newcommand{\bi}{\begin{exercise}}
\newcommand{\ei}{\end{exercise}}
\newcommand{\bo}{\begin{proposition}}
\newcommand{\eo}{\end{proposition}}
\newcommand{\br}{\begin{remark}}
\newcommand{\er}{\end{remark}}
\newcommand{\be}{\begin{equation}}
\newcommand{\ee}{\end{equation}}
\newcommand{\ba}{\begin{align}}
\newcommand{\ea}{\end{align}}
\newcommand{\bn}{\begin{enumerate}}
\newcommand{\en}{\end{enumerate}}
\newcommand{\bg}{\begin{align*}}
\newcommand{\bcs}{\begin{cases}}
\newcommand{\ecs}{\end{cases}}

\newcommand{\NN}{{\mathbb N}}

\newcommand{\bean}{\begin{eqnarray*}}
\newcommand{\eean}{\end{eqnarray*}}

\makeatletter % @ is now a normal "letter" for TeX
\renewcommand\theequation{\thesection.\arabic{equation}}
\@addtoreset{equation}{section}
\makeatother%@isrestoredasa"non-letter"characterforTeX
\numberwithin{equation}{section}

\begin{document}

\begin{center}
\textbf{Infinitely many solutions for  elliptic system with Hamiltonian type}\\
\end{center}

%\begin{center}
%Weimin Zhang\\
%%$^{1}$ School of Mathematics and Computer Science, Shanxi Normal University,\\ Taiyuan, 030031, P.R. China\\
%School of Mathematical Sciences,  Key Laboratory of Mathematics and Engineering Applications (Ministry of Education) \& Shanghai Key Laboratory of PMMP,  East China Normal University, Shanghai 200241, China
%\end{center}

\begin{center}
Jia Zhang$^1$,  Weimin Zhang$^{2, *}$\\
\medskip

$^1$ School of Mathematical Sciences,  Key Laboratory of Mathematics and Engineering Applications (Ministry of Education) \& Shanghai Key Laboratory of PMMP,  East China Normal University, Shanghai, 200241, P.R. China
\smallskip

$^2$ School of Mathematical Sciences, Zhejiang Normal University, Jinhua, 321004, P.R. China
\end{center}

\begin{center}
\renewcommand{\theequation}{\arabic{section}.\arabic{equation}}
\numberwithin{equation}{section}
\footnote[0]{\hspace*{-7.4mm}
$^*$Corresponding author.
{E-mail addresses: 52265500023@stu.ecnu.edu.cn (J. Zhang), zhangweimin2021@gmail.com (W. Zhang).}}
\end{center}
%\address[W.~M.~Zhang]{\newline\indent School of Mathematical Sciences
%\newline\indent
%East China Normal University
%\newline\indent
%Shanghai 200241, P.R. China}
%\email{\href{mailto:52205500026@stu.ecnu.edu.cn}{52205500026@stu.ecnu.edu.cn}}
%\address[H.~R.~Sun]{\newline\indent School of Mathematics and Statistics
%\newline\indent
%Lanzhou University,
%\newline\indent Lanzhou 730000, P.R. China}
%\email{\href{mailto:hrsun@lzu.edu.cn}{hrsun@lzu.edu.cn}}
%\address[J.~Zhang]{\newline\indent College of Mathematics and Statistics
%\newline\indent
%Chongqing Jiaotong University
%\newline\indent
%Chongqing 400074, P.R. China}
%\email{\href{mailto:zhangjianjun09@tsinghua.org.cn}{zhangjianjun09@tsinghua.org.cn}}
%\address[Z.~F.~Jin]{\newline\indent School of Mathematics and Computer Science
%\newline\indent
%Shanxi Normal University,
%\newline\indent Taiyuan, Shanxi 030031, China}
%\email{\href{mailto:jinzhf15@lzu.edu.cn}{jinzhf15@lzu.edu.cn}}
%\thanks{(1) Corresponding author: \texttt{hrsun@lzu.edu.cn}}
%\thanks{(2) H.~R.~Sun was partly supported by the NSFC (Grants No. 11671181) and NSF of Gansu Province of China (Grants No. 21JR7RA535)}
%\thanks{(3) J.~J.~Zhang was supported by NSFC(No.11871123)}
%\subjclass[2010]{35A15, 35J60, 58E05}
\begin{abstract}
In this paper, we use Legendre-Fenchel transform and a space decomposition to carry out Fountain theorem and dual Fountain theorem for the following elliptic system of Hamiltonian type:
\[
\begin{cases}
\begin{aligned}
-\Delta u&=H_v(u, v) \,\quad&&\text{in}~\Omega,\\
-\Delta v&=H_u(u, v) \,\quad&&\text{in}~\Omega,\\
u,\,v&=0~~&&\text{on} ~ \partial\Omega,\\
\end{aligned}
\end{cases}
\]
where $N\ge 1$, $\Omega \subset \mathbb{R}^N$ is a bounded domain and $H\in C^1( \mathbb{R}^2)$ is strictly convex, even and subcritical. We mainly present two results: (i) When $H$ is superlinear, the system has infinitely many solutions, whose energies tend to infinity. (ii) When $H$ is sublinear, the system has infinitely many solutions, whose energies are negative and tend to 0. As a byproduct, the Lane-Emden system under subcritical growth has infinitely many solutions.
\end{abstract}


\section{Introduction}\label{s1}
 In this paper, we are concerned with the existence of infinitely many solutions   for the following elliptic system of Hamiltonian type:
\begin{equation}\label{2305191141}
\begin{cases}
\begin{aligned}
-\Delta u&=H_v(u, v) \,\quad&&\text{in}~\Omega,\\
-\Delta v&=H_u(u, v) \,\quad&&\text{in}~\Omega,\\
u,\,v&=0~~&&\text{on} ~ \partial\Omega,\\
\end{aligned}
\end{cases}
\end{equation}
where $N\ge 1$, $H\in C^1(\mathbb{R}^2)$ and $\Omega\subset \mathbb{R}^N$ is a bounded domain. Formally, the associated energy functional of \eqref{2305191141} is 
\[
\mathcal{I}(u, v)=\int_{\Omega}\nabla u \nabla v dx-\int_{\Omega}H(u, v) dx.
\]
It is noteworthy that the quadratic part of $\mathcal{I}$ is strongly indefinite, namely, it is neither positive definite nor negative definite on usual energy spaces. The strong indefiniteness prevents people from applying mountain pass theorem, creating many difficulties while dealing with the existence issue via the variational method. Various methods have been developed to overcome such difficulties. Benci-Rabinowitz \cite{BR1979} established a linking theorem in Hilbert framework for strongly indefinite functionals. With this result, when the functional $\mathcal{I}$ is defined on $H_0^1(\Omega)\times H_0^1(\Omega)$, the minimax method is feasible, but this space imposes great restrictions for the growth of $H$. Previous works often concern whether solutions of \eqref{2305191141} exist under the growth
\begin{equation}\label{24103016244}
|H(u, v)|\le C(|u|^{p+1}+|v|^{q+1})
\end{equation}
with
\begin{equation}\label{2410301620}
\frac1{p+1}+\frac1{q+1}\ge \frac{N-2}{N},\quad p, q>0.
\end{equation}
 In this case, the energy space could be chosen as  $W_0^{1, t}(\Omega)\times W_0^{1, \frac{t}{t-1}}(\Omega)$ for some $t>1$. However, Benci-Rabinowitz's linking theorem could not work in a general Banach setting, which makes such problem more challenging. Recall that when $N\ge 3$, the curve
\begin{equation}\label{2501091719}
\frac1{p+1}+\frac1{q+1}= \frac{N-2}{N}
\end{equation}
is called the critical hyperbola. When the inequality \eqref{2410301620} is strict, we call $H$ possesses subcritical growth. In fact, \eqref{2501091719} is a borderline for the existence of problem \eqref{2305191141}, since a Rellich identity in \cite{Mitidieri1993} can yield a non-existence result for critical and supercritical case.

\medskip
There are many articles devoted to finding nontrivial solutions of \eqref{2305191141} under subcritical or critical conditions.  For the subcritical case, Hulshof-van der Vorst \cite{Hv1993} used the interpolation space framework to deal with the system with $H(u, v) = F(u)+G(v)$. de Figueiredo-Felmer \cite{dF1994} also used the interpolation space method to deal with the general Hamiltonian system, but the $p, q$  have the restriction:
\begin{equation}
(N-4)\max\{p, q\}<N+4.
\end{equation}
 Cl\'ement-van der Vorst in \cite{CV1995} proposed dual method to study  \eqref{2305191141}, which is based on Legendre-Fenchel transform. This method can transform the strongly indefinite problem into a problem with mountain pass structure, but the strict convexity of  $H$ is required. de Figueiredo-do \'{O}-Ruf \cite{ddR_JFA2005} solved \eqref{2305191141} in finite-dimensional subspaces of a Sobolev-Orlicz space, by an approximation, they found a nontrivial solution of \eqref{2305191141}. Additionally, Hulshof-Mitidieri-van der Vorst \cite{HMv1998} used the dual method to treat the critical case.
\medskip

Very few papers discuss the multiplicity of solutions of \eqref{2305191141}. It is well-known that Fountain theorem and dual Fountain theorem established respectively by Bartsch \cite{B1993} and Bartsch-Willem \cite{BW1995} are often used to find infinitely many solutions, see for instance \cite{Bd1999, Liu2010, FZ2003}. Bartsch-de Figueiredo \cite{Bd1999} applied Fountain theorem to the Hamiltonian system with interpolation space setting, whose functional is defined on a Hilbert space. Their models have some restrictions on Hamiltonian systems, for instance, that cannot be simply applied to all Lane-Emden systems with subcritical growth.  In this paper, by means of Legendre-Fenchel transform and a space decomposition method, we will carry out Fountain theorem and dual Fountain theorem to get infinitely many solutions of \eqref{2305191141}. As a byproduct, the Lane-Emden system 
\begin{equation}\label{2401050900}
\begin{cases}
\begin{aligned}
-\Delta u&=|v|^{q-1}v \,\quad&&\text{in}~\Omega,\\
-\Delta v&=|u|^{p-1}u \,\quad&&\text{in}~\Omega,\\
u,\,v&=0~~&&\text{on} ~ \partial\Omega,\\
\end{aligned}
\end{cases}
\end{equation}
under subcritical growth ($pq\neq 1$) has infinitely many solutions.

\subsection{Main results}

\medskip
In this paper, we will impose the following assumptions on $H$:
\begin{itemize}
\item[$(H1)$] There exist $p,\, q>0$ and  positive numbers $C_1, C_2$ such that 
\[
\begin{split}
&C_1|u|^{p+1}\le H_u(u,v)u\le C_2\left( |u|^{p+1}+  |u|^{\alpha}|v|^{\beta}\right),\\
&C_1|v|^{q+1}\le H_v(u,v)v\le C_2\left( |v|^{q+1}+  |u|^{\alpha}|v|^{\beta}\right),
\end{split}
\]
 with
\begin{equation}\label{2402051803}
\frac{\alpha}{p+1}+\frac{\beta}{q+1}=1,\;\; \alpha, \beta>1;
\end{equation}
\item[$(H2)$] $H$ is strictly convex;
\item[$(H3)$] There are $\theta \in (0, 1)$ and positive numbers $C_3$, $C_4$ such that
\[
\theta H_u(u,v)u+ (1- \theta)H_v(u,v)v- H(u,v)\ge C_3\Big(|u|^{p+1}+|v|^{q+1}\Big)-C_4;
\]
\item[$(H4)$] $H(-u, -v)=H(u, v)$ for any $(u, v)\in \mathbb{R}^2$.
\end{itemize}

\begin{remark}
We emphasize that $(H1)$ and $(H2)$ have some restrictions on $H$. However, there are still  many examples satisfying such conditions, for instance
\begin{equation}\label{2402052125}
H_\varepsilon(u, v)=|u|^{p+1}+|v|^{q+1}+\varepsilon |u|^{\alpha}|v|^{\beta}
\end{equation}
with $|\varepsilon|$ small enough and $\alpha,\, \beta$ satisfy \eqref{2402051803}. In particular, the Lane-Emden system \eqref{2401050900} satisfies $(H1)$-$(H4)$ if we choose $\varepsilon=0$.
\end{remark}

We proceed by dual method. Set the function space $X:=L^{1+\frac1p}(\Omega)\times L^{1+\frac1q}(\Omega)$ and its dual space is $X^*=L^{p+1}(\Omega)\times L^{q+1}(\Omega)$. We define the functional
\begin{equation}\label{2011011909}
\mathcal{J}(f, g)=\mathcal{H}^*(f, g)-\int_{\Omega}g\mathcal{A}f dx,\quad \forall \, (f, g)\in X,
\end{equation}
where $\mathcal{A}$ is the inverse of $-\Delta: W_0^{1, 1+\frac1p}(\Omega)\cap W^{2, 1+\frac1p}(\Omega)\to L^{1+\frac1p}(\Omega)$ and $\mathcal{H}^*$ denotes the Legendre-Fenchel transform of 
\begin{equation}\label{2411011908}
\mathcal{H}(u, v)=\int_{\Omega}H(u, v) dx, \quad \forall\; (u, v)\in X^*,
\end{equation}
seeing subsection \ref{2501052108} below. Note that \eqref{2011011909} and \eqref{2411011908} are well defined due to $(H1)$, \eqref{2410301620} and $\mathcal{A}f\in W_0^{1, 1+\frac1p}(\Omega)\cap W^{2, 1+\frac1p}(\Omega)$ for any $f\in L^{1+\frac1p}$. $\mathcal{A}$ is also well defined in $L^r(\Omega)$ for all $r\ge 1$. Without confusion, we will not explain the domain of $\mathcal{A}$. 

\smallskip
With the assumptions $(H1)$-$(H2)$, the functional $\mathcal J$ is differentiable and any critical point of $\mathcal J$ corresponds to a critical point of $\mathcal I$, hence it corresponds to a solution of \eqref{2305191141}, seeing subsection \ref{2501052108}. We will prove in Lemma \ref{2401052112} that $\mathcal{J}(f, g)$ is a critical value of $\mathcal I$ for any critical point $(f, g)$ of $\mathcal J$. 

\smallskip

Under the superlinear and subcritical condition on $H$, i.e.,  $p,\, q$ in $(H1)$-$(H4)$ satisfy
\begin{equation}\label{2410222035}
1>\frac1{p+1}+\frac1{q+1}>\frac{N-2}{N},
\end{equation}
we will first apply Fountain theorem to obtain infinitely many solutions of \eqref{2305191141}.
Since $\mathcal J$ has two asymmetric variables, we encounter some difficulties while constructing the linking structure. On the one hand, we need to establish a space decomposition for the space $X$, which is used to construct linking geometry. By means of the method in \cite{YZ_ANS2024}, we obtain a subtle decomposition for Banach space with two asymmetric variables in subsection \ref{2501091039}. On the other hand, another important part of linking structure is the intersection property. Observing that the asymmetric variables have a great effect on the minimax process, our linking sets actually are not located on a finite-dimensional linear space. For this reason, the classical topological degree could not deduce directly the  intersection property. For that, we will introduce the notion of $\mathbb{Z}_2$-cohomological index in subsection \ref{2501091050}, whose fine properties conclude the desired result in Proposition \ref{2410121928}.

\medskip
Our first main result can be stated as follows.
\begin{theorem}\label{thm1}
Assume that $H\in C^1(\mathbb{R}^2)$ satisfies $(H1)$-$(H4)$ and $p,\, q>0$ satisfy \eqref{2410222035}. Then there exist infinitely many nontrivial solutions $(u_j, v_j)$ to \eqref{2305191141} such that $\mathcal{I}(u_j, v_j)\to \infty$ as $j\to\infty$.
\end{theorem}

Next, we are concerned with the coercive case, i.e. $p,\, q$ in $(H1)$, $(H2)$ and $(H4)$ satisfy
\begin{equation}\label{2410222037}
\frac1{p+1}+\frac1{q+1}>1.
\end{equation}
This case is different from the case of \eqref{2410222035}. We will carry out the steps of dual Fountain theorem for the functional $\mathcal{J}$ in the space $X$. As we explained previously, while establishing the linking geometry and intersection property, our main difficulty is still that the energy functional has two asymmetric variables.

\medskip
The second main result is the following.
\begin{theorem}\label{thm2}
Assume that $H\in C^1(\mathbb{R}^2)$ satisfies $(H1)$, $(H2)$, $(H4)$ and $p,\,q>0$ satisfy \eqref{2410222037}. Then there exist infinitely many nontrivial solutions $(u_j, v_j)$ to \eqref{2305191141} such that $\mathcal{I}(u_j, v_j)\to 0^-$ as $j\to\infty$.
\end{theorem}

\begin{remark}
The assumption $(H3)$ is an important condition in Theorem \ref{thm1}, which could guarantee the boundedness of Palais-Smale sequence. However, we do not need this assumption in the coercive case, since the minimax sequence is naturally bounded.
\end{remark}
\begin{remark}
In Theorem \ref{thm2}, since the energy functional $\mathcal{J}$ is coercive, it is still valid to use $\mathbb{Z}_2$-genus to construct minimax. If so, we can get infinitely many solutions due to Lusternik-Schnirelman theory, but it may not conclude the energy tendency. 
\end{remark}

\medskip

The rest of the paper is organized as follows. In Section \ref{2501091922}, we introduce some preliminary results, including the notion of Legendre-Fenchel transform, $\mathbb{Z}_2$-cohomological index and space decomposition. In Section \ref{2501091925}, we are dedicated to proving Theorem \ref{thm1}. In Section \ref{2411022316}, we are devoted to showing Theorem \ref{thm2}.

\section{Preliminary}\label{2501091922}
Throughout this paper, $C, C_1, C_2,...$ always denote generic positive constants  and $|\cdot|_r$ means  the usual norm of $L^r(\Omega)$ for $r \geq 1$.
\subsection{Legendre-Fenchel transform}\label{2501052108}
Let $V$ be a Banach space. For a function $G: V\to \mathbb{R}\cup +\infty$, $G\not\equiv +\infty$, its Legendre-Fenchel transform is given by
\begin{equation}\label{2305291511}
G^*(u^*)=\sup\{\langle u^*, v \rangle-G(v): v\in V\},\quad \forall\, u^*\in V^*,
\end{equation}
where $V^*$ is the dual space of $V$, and $\langle \cdot, \cdot \rangle$ is the dual operation. 
\begin{lemma}\label{2410222032}
Assume that $G$ is an even functional on a Banach space $V$. Then $G^*$ is also even.
\end{lemma}
\begin{proof}
According to the definition \eqref{2305291511} and the assumption that $G$ is even,
\[
\begin{aligned}
G^*(-u^*)&=\sup\{\langle -u^*, v \rangle-G(v): v\in V\}\\
&=\sup\{\langle u^*, -v \rangle-G(-v): v\in V\}\\
&=\sup\{\langle u^*, v \rangle-G(v): v\in V\}\\
&=G^*(u^*).
\end{aligned}
\]
The proof is completed.
\end{proof}
%\begin{lemma}\label{2410272114}
%Let $G_1$, $G_2$ be two functionals on a Banach space $V$, then for any $\theta\in (0, 1)$,
%$$(G_1+G_2)^*(u^*)\le G_1^*(\theta u^*)+G_2^*((1-\theta)u^*),\quad \forall\, u^*\in V^*.$$
%\end{lemma}
%\begin{proof}
%Set $G=G_1+G_2$. According to the definition of Legendre-Fenchel transform, for any $u^*\in V^*$,
%\begin{equation}\label{2410272029}
%\begin{aligned}
%G^*(u^*)=&\sup\{\langle u^*, v \rangle-G_1(v)-G_2(v): v\in V\}\\
%\le &\sup\{\langle \theta u^*, v \rangle-G_1(v)\}+\sup\{\langle (1-\theta)u^*, v \rangle-G_2(v): v\in V\}\\
%=& G_1^*(\theta u^*)+G_2^*((1-\theta)u^*).
%\end{aligned}
%\end{equation}
%Hence we complete the proof.
%\end{proof}

We denote by $H^*$ the Legendre-Fenchel transform of $H:\mathbb{R}^2\to \mathbb{R}$, that is
\[
H^*(u, v)=\underset{(t, s)\in \mathbb{R}^2}{\sup}\left\{ tu+sv-H(t, s)\right\},\quad \forall\; (u,v)\in \mathbb{R}^2,
\]
and $\mathcal{H}^*$ denotes the Legendre-Fenchel transform of $\mathcal{H}: X^*\to \mathbb{R}$, namely
\begin{equation}
\mathcal{H}^*(f, g)=\underset{(u, v)\in X^*}{\sup} \left\{ \int_{\Omega} (fu+gv) dx-\mathcal{H}(u, v)\right\},\quad \forall\; (f, g)\in X.
\end{equation}
 Next, the lemmas below can refer to \cite[Section 4]{Z_arXiv2024}.

\begin{lemma}[\cite{Z_arXiv2024} Lemma 4.2 ]\label{2211011024}
Assume that $H$ is convex and satisfies $(H1)$, then there are positive constants $A_1$, $A_2$ such that
\[
A_1 \Big(|f|_{1+\frac1p}^{1+\frac1p}+ |g|_{1+\frac1q}^{1+\frac1q}\Big)\le \mathcal{H}^*(f, g)\le A_2 \Big(|f|_{1+\frac1p}^{1+\frac1p}+ |g|_{1+\frac1q}^{1+\frac1q}\Big), \quad \forall\; (f, g) \in X.
\]
\end{lemma}
\begin{lemma}[\cite{Z_arXiv2024} Lemmas 4.4, 4.5 ]\label{22112307}
Assume that $H$ satisfies $(H1)$-$(H2)$, then $H^*\in C^1(\mathbb{R}^2, \mathbb{R})$ and $\mathcal{H}^{*}\in C^1(X, \mathbb{R})$. More precisely,
\begin{equation}\label{2211052325}
\mathcal{H}^{*}(f, g)=\int_{\Omega}H^*(f, g) dx, \quad 
\left\langle (\mathcal{H}^{*})'(f, g), (\widetilde{f}, \widetilde{g}) \right\rangle=\int_{\Omega}H^*_u(f, g) \widetilde{f}+H^*_v(f, g)\widetilde{g} dx.
\end{equation}
\end{lemma}

\begin{lemma}\label{2410221530}
Under the conditions $(H1)$-$(H2)$, if $(f, g)\in X$ is a critical point of $\mathcal{J}$, then $\nabla H^*(f, g)$ is a classical solution of \eqref{2305191141}.
\end{lemma}
\begin{proof}
A similar discussion can be found in \cite{CV1995} and \cite[Section 4]{Z_arXiv2024}. If $(u, v)\in X$ is a weak solution of \eqref{2305191141}, when $H$ has subcritical growth, then $(u, v)\in L^{\infty}(\Omega)\times L^{\infty}(\Omega)$, seeing for example \cite[Lemma 3.9]{Z_arXiv2024}, hence that are smooth.
\end{proof}

\begin{lemma}\label{2401052112}
Under the condition $(H1)$-$(H2)$, we sssume that  $(f, g)$ is a critical point of $\mathcal{J}$, and let $(u, v)=\nabla H^*(f, g)$, then
\begin{equation}\label{2410221548}
\mathcal{I}(u, v)=\mathcal{J}(f, g).
\end{equation}
\end{lemma}
\begin{proof}
Since $(u, v)=\nabla H^*(f, g)$,  we follow from \cite[Chapter I, Lemma 6.3]{Struwe2008} that $(f, g)=\nabla H(u, v)$. Hence 
\begin{equation}\label{2410221541}
\int_{\Omega}f\mathcal{A}g dx=\int_{\Omega}\nabla u \nabla v dx.
\end{equation}
Using \cite[Chapter I, Lemma 6.3]{Struwe2008} again, that is
\[
\int_{\Omega}\left(fu +gv\right) dx=\int_{\Omega}H(u, v) + H^*(f, g)\, dx,
\] 
then we get
\begin{equation}\label{2410221542}
\begin{aligned}
\int_{\Omega} H^*(f, g) dx&=\int_{\Omega}H_u(u, v)u dx+\int_{\Omega}H_v(u, v)v dx-\int_{\Omega}H(u, v) dx\\
&=2\int_{\Omega}\nabla u \nabla v dx-\int_{\Omega}H(u, v) dx.
\end{aligned}
\end{equation}
By \eqref{2410221541} and \eqref{2410221542}, we can obtain \eqref{2410221548}. 
\end{proof}

\subsection{$\mathbb{Z}_2$-cohomological index}\label{2501091050}
To establish the linking structure, we shall prove some intersection properties in Propositions \ref{2410121928} and \ref{2410132041}, where we will use some properties of $\mathbb{Z}_2$-cohomological index. Recall the definition of $\mathbb{Z}_2$-cohomological index by Fadell and Rabinowitz \cite{FR1978}. For a Banach space $W$, let $M$ be a symmetric subset of $W\backslash\{0\}$. Let $\overline{M}=M/\mathbb{Z}_2$ be the quotient space (in which $u$ and $-u$ are identified), and let $\phi: \overline{M}\to \mathbb{R}P^{\infty}$ be the classifying map of the space $\overline{M}$, which induces a homomorphism $\phi^*: H^*(\mathbb{R} P^\infty)\to H^*(\overline{M})$ of Alexander-Spanier cohomology rings with coefficients in $\mathbb{Z}_2$. We may identify $H^*(\mathbb{R} P^\infty)$ with the polynomial ring $\mathbb{Z}_2[\omega]$. The cohomological index of $M$ is 
\[
{\rm Ind}(M)=
\begin{cases}
 \sup\{k\in \NN: \phi^*(\omega^k)\neq 0\}  \quad \text{if}~M\neq \varnothing,\\
0~~~~~~~~~~~~~~~~~~~~~~~~~~~~~\,~~~~\text{if}~M= \varnothing.
\end{cases}
\]
We denote $\mathcal{M}$ be the class of symmetric subsets of $W\backslash\{0\}$. The cohomological index is a map ${\rm Ind}:\mathcal{M}\to \mathbb{N}\cup \{0, \infty\}$ and possesses the following properties.

\begin{lemma}[\cite {PSY2016} Proposition 2.1]\label{2410222003}
\begin{itemize}
\item[\rm (i)] ${\rm Ind}(A)=0$ if and only if $A=\varnothing$.
\item[\rm (ii)] If there exists an odd map from $A$ to $B$, then ${\rm Ind}(A)\le {\rm Ind}(B)$. In particular, ${\rm Ind}(A)\le {\rm Ind}(B)$ if $A\subset B$.
\item[\rm (iii)] ${\rm Ind}(A)\le {\rm dim}(W)$.
\item[\rm (iv)]  If $A$, $A_0$ and $A_1$ are closed, and $\varphi:A\times [0, 1]\to A_0\cup A_1$ is a continuous map such that $\varphi(-u, t)=-\varphi(u, t)$ for all $(u, t)\in A\times [0, 1]$, $\varphi(A\times [0, 1])$ is closed, $\varphi(A\times\{0\})\subset A_0$ and $\varphi(A\times\{1\})\subset A_1$, then
\[
{\rm Ind}({\rm Im}(\varphi)\cap A_0\cap A_1)\ge {\rm Ind}(A).
\]
\item[\rm (v)] If $U$ is a bounded closed symmetric neighborhood of $0$, then ${\rm Ind}(\partial U)={\rm dim}(W)$.
\end{itemize}
\end{lemma}

%\subsection{Palais-Smale condition and deformation lemma}
%\begin{definition}
% Let $I\in C^1(X, \mathbb{R})$ and $c\in \mathbb{R}$. The functional $I$ satisfies the $(PS)_c$ condition if any sequence $\{u_{j}\}\subset X$ such that 
%\begin{equation}
%I(u_{j})\to c,\; I'(u_{j})\to 0\; \mbox{as } j\to \infty
%\end{equation}
%admits a subsequence converging to a critical point of $I$.
% \end{definition}

\subsection{Space decomposition}\label{2501091039}
Another important ingredient for establishing the linking structure is to find some suitable space decompositions. Due to the fact that our functional possesses two asymmetric variables, the space decomposition is more complex in applying Fountain theorem and its dual version.
\medskip

Applying the decomposition result \cite[Theorem 5.2]{YZ_ANS2024} to $L^{1+\frac1q}(\Omega)$ and its dual space $L^{q+1}(\Omega)$, there exist two sequences $\{v_j\}\subset L^{q+1}(\Omega)$ and $\{g_j\}\subset L^{1+\frac1q}(\Omega)$ such that
\[
L^{q+1}(\Omega)=\overline{{\rm span}\{v_j, j\ge 1\}}
\]
and
 \begin{equation}\label{2308271605}
\forall\; j\ge 1,\quad \int_{\Omega}g_j v_j dx=1,\quad \int_{\Omega}g_j v_k dx=0 \;\;\mbox{for all }k\neq j.
\end{equation}
 Furthermore, if we let
\[
E_j = {\rm span}\{g_k, 1\le k \le j\},\;\;E_j^{\perp}=\underset{1 \le k\le j}{\cap}\left\{g\in L^{1+\frac1q}(\Omega): \int_{\Omega}g v_k dx=0\right\}, \quad \forall\; j \ge 1,
\]
then 
\begin{itemize}
\item $L^{1+\frac1q}(\Omega)= E_j \oplus E_j^{\perp}=\overline{\cup_{j\ge 1} E_j}$, ${\rm dim}(E_j)=j$;
\item Let $\{w_j\} \subset L^{1+\frac1q}(\Omega)$ be a bounded sequence such that $w_j\in E_j^{\perp}$ for any $j\ge 1$, then $w_j\to 0$ weakly in $L^{1+\frac1q}(\Omega)$;
\item We denote by $P_j$ the projection from $L^{1+\frac1q}(\Omega)$ into $E_j$ parallel to $E_j^\perp$, then for any $\varphi\in L^{1+\frac1q}(\Omega)$, it holds that $P_j \varphi\to\varphi$ in $L^{1+\frac1q}(\Omega)$.
\end{itemize}

\begin{lemma}\label{2411012206}
 $\{\mathcal{A}g_j\}$ is a basis of $L^{p+1}(\Omega)$. In other words, $\{\mathcal{A}g_j\}$ is linearly independent and $L^{p+1}(\Omega)=\overline{{\rm span}\{\mathcal{A}g_j, j\ge 1\}}$.
\end{lemma}
\begin{proof}
Since $\{g_j\}$ is linearly independent in $L^{q+1}(\Omega)$, so does $\{\mathcal{A}g_j\}$ in $W_0^{1, 1+\frac1q}(\Omega)\cap W^{2, 1+\frac1q}(\Omega)$ and $L^{p+1}(\Omega)$.  By density property, it remains to prove $C_c^{\infty}(\Omega)\subset \overline{{\rm span}\{\mathcal{A}g_j, j\ge 1\}}$. For any $u\in C_c^{\infty}(\Omega)$, obviously $g:=-\Delta u\in L^{1+\frac1q}(\Omega)$. In view of the properties of the projection operator $P_j$, if we let $w_j:=P_jg$, there holds that
\[
w_j\to g\quad\mbox{in }L^{1+\frac1q}(\Omega).
\]
Hence
\[
\mathcal{A}w_j\to \mathcal{A}g=u\quad\mbox{in }L^{p+1}(\Omega).
\]
Since $\mathcal{A}w_j\in {\rm span}\{\mathcal{A}g_k, 1\le k\le j\}$, the proof is done.
\end{proof}

As the proof of \cite[Theorem 5.2]{YZ_ANS2024}, by Lemma \ref{2411012206}, we can easily find a basis $\{f_j\}$ of $L^{1+\frac1p}(\Omega)$ such that 
\begin{equation}\label{2408291506}
\forall\; j\ge 1,\quad \int_{\Omega}f_j\mathcal{A}g_j dx=1,\quad \int_{\Omega}f_j\mathcal{A}g_k dx=0 \;\;\;\;\mbox{for any }k\neq j.
\end{equation}
Similar to the processes above, we denote
\[
\widetilde{E}_j = {\rm span}\{f_k, 1\le k \le j\},\;\;\widetilde{E}_j^{\perp}=\underset{1 \le k\le j}{\cap}\left\{f\in L^{1+\frac1p}(\Omega): \int_{\Omega}f\mathcal{A} g_k dx=0\right\},  \quad \forall\; j \ge 1.
\]
Thus 
\[
 L^{1+\frac1p}(\Omega)= \widetilde{E}_j \oplus \widetilde{E}_j^{\perp},\quad {\rm dim}(\widetilde{E}_j)=j\quad \mbox{ and } \quad L^{1+\frac1p}(\Omega)=\overline{\cup_{j\ge 1} \widetilde{E}_j}.
\]
We denote by $ \widetilde{P}_j$ the projection from $L^{1+\frac1p}(\Omega)$ into $\widetilde{E}_j$ parallel to $\widetilde{E}_j^\perp$. Also, $\widetilde{P}_j \varphi\to\varphi$ in $L^{1+\frac1p}(\Omega)$ as $j\to\infty$ for any $\varphi\in L^{1+\frac1p}(\Omega)$.

\medskip
Now let us define two sequences 
\begin{equation}\label{2412220128}
\alpha_j:=\underset{g\in E_j^{\perp},\, |g|_{1+\frac1q}=1}{\sup}|\mathcal{A}g|_{p+1},\quad \beta_j:=\underset{f\in \widetilde{E}_j^{\perp},\, |f|_{1+\frac1p}=1}{\sup}|\mathcal{A}f|_{q+1}.
\end{equation}
\begin{lemma}\label{2408292013}
$\alpha_j,\,\beta_j$ are finite and $\alpha_j,\,\beta_j\to 0$ as $j\to\infty$.
\end{lemma}
\begin{proof}
Since $p,\, q$ are subcritical and using Sobolev embedding, $\alpha_j$ and $\beta_j$ are finite. Next, we prove $\alpha_j\to 0$. By contradiction, we assume that there exist $w_j\in E_j^{\perp}$ and $c>0$ such that
\begin{equation}\label{2408221043}
|w_j|_{1+\frac1q}=1,\quad |\mathcal{A}w_j|_{p+1}>c.
\end{equation}
As we stated before, $w_j$ weakly converges to 0 in $L^{1+\frac1q}(\Omega)$. Since the  map $\mathcal{A}: L^{1+\frac1q}(\Omega)\to L^{p+1}(\Omega)$ is compact, $|\mathcal{A}w_j|_{p+1}\to 0$, which contradicts \eqref{2408221043}. So $\alpha_j\to 0$, and similarly $\beta_j\to 0$.
\end{proof}

We now define for $n\ge 1$,
\begin{equation}\label{2408291509}
F_{n}:={\rm span}\{(f_j, g_j), 1\le j\le n\},
\end{equation}
\[
G_n:=\widetilde{E}_n\times E_n,\quad G_{n}^{\perp}:=\widetilde{E}_n^{\perp}\times E_n^{\perp}.
\]
For $m\ge n\ge 1$, we denote 
\[
\widetilde{E}_n^m:={\rm span}\{f_j: n\le j\le m\},\quad  E_n^m:={\rm span}\{g_j: n\le j\le m\}.
\]
and
\[
G_n^m:=\widetilde{E}_n^m\times E_n^m.
\]

\begin{lemma}\label{2410140905}
For any $m\ge n\ge 1$, there hold that ${\rm dim}(F_n)=n$, ${\rm dim}(G_n)=2n$, ${\rm dim}(G_n^m)=2m-2n+2$ and $X=G_n\oplus G_n^{\perp}$.
\end{lemma}
\begin{proof}
It is obvious that ${\rm dim}(F_n)=n$. Next, one can verify that
\[
\{(f_j, 0),\; (0, g_j): \;\; n\le j\le m\}
\]
is a basis of $G_n^m$. Thus, ${\rm dim}(G_n)=2n$ and ${\rm dim}(G_n^m)=2m-2n+2$. 
 
 \smallskip
 Finally, we prove the last claim. For any $(f, g)\in L^{1+\frac1p}(\Omega)\times L^{1+\frac1q}(\Omega)$, we have
\[
(f, g)=(\widetilde{P}_n f, P_n g)+(f-\widetilde{P}_n f, g-P_ng).
\]
So $L^{1+\frac1p}(\Omega)\times L^{1+\frac1q}(\Omega)=G_n+G_n^{\perp}$. To prove that the plus at $G_n+G_n^{\perp}$ is a direct sum, it remains to verify $G_n\cap G_n^{\perp}=\{0\}$, which is obvious. 
\end{proof}

\begin{lemma}\label{2408292319}
For any $n\in \mathbb{N}^+$, there holds that
\[
\left|\int_{\Omega}f\mathcal{A}g dx\right|\le \gamma_n|f|_{1+\frac1p}|g|_{1+\frac1q},\quad \forall \, (f, g)\in G_n^{\perp},
\]
where $\gamma_n=\min\{\alpha_n, \beta_n\}$.
\end{lemma}
\begin{proof}
For any $(f, g)\in G_n^{\perp}$, we have
\[
\begin{aligned}
\left|\int_{\Omega}f\mathcal{A}g dx\right|\le |f|_{1+\frac{1}{p}}|\mathcal{A}g|_{1+p}\le \alpha_n|f|_{1+\frac{1}{p}}|g|_{1+\frac1q}.
\end{aligned}
\]
On the other hand, since 
\[
\int_{\Omega}f\mathcal{A}g dx=\int_{\Omega}\mathcal{A}fg dx,
\]
we obtain
\[
\begin{aligned}
\left|\int_{\Omega}f\mathcal{A}g dx\right|\le |\mathcal{A}f|_{1+q}|g|_{1+\frac1q}\le \beta_n|f|_{1+\frac{1}{p}}|g|_{1+\frac1q}.
\end{aligned}
\]
The above two inequalities can deduce the conclusion.
\end{proof}
\begin{lemma}\label{2408300019}
For any $n\in \mathbb{N}^+$, there exists $C(n)>0$ such that
\[
\int_{\Omega}f\mathcal{A}g dx\ge C(n), \quad \forall \, (f, g)\in F_n\mbox{ with }\|(f, g)\|_{X}=1.
\]
\end{lemma}
\begin{proof}
If otherwise, we will get a minimizing sequence $\{(w_j, y_j)\}\subset F_n$ such that as $j\to\infty$,
\[
\int_{\Omega}w_j\mathcal{A}y_j dx\to 0.
\]
Since $F_n$ is a finite dimensional space and $\{(w_j, y_j)\}$ is a bounded sequence, there exists some $(w, y)\in F_n$ such that
\[
\int_{\Omega}w\mathcal{A}y dx= 0,\quad \mbox{and}\quad \|(w, y)\|_{X}=1.
\]
This is a contradiction since $\mathcal{A}$ is positive definite in $F_n$.
\end{proof}
\section{Superlinear case: infinitely many solutions with positive energies}\label{2501091925}
In this section, we mainly use Fountain theorem to find infinitely many solutions of \eqref{2305191141} with energies tending to infinity. Thoughout this section, we assume that $p,\, q$ are superlinear and subcritical, i.e. $p,\, q>0$ satisfy \eqref{2410222035}. Thus, there exist $k, l>1$ satisfying
\begin{equation}\label{2211012130}
\frac{p}{p+1}>\frac{k}{k+l}, \quad \frac{q}{q+1}>\frac{l}{k+l}.
\end{equation}
We  denote
\[
\begin{aligned}
&B^n_\rho:=\{(t^{k}f, t^{l}g): (f, g)\in F_{n}, \,\|(f, g)\|_X= 1,\; 0\le t\le \rho\},\\
&Q^n_{\rho}:=\{(\rho^{k}f, \rho^{l}g): (f, g)\in F_{n}, \,\|(f, g)\|_X=1\},
\end{aligned}
\]
and
$$S^n_\rho:=\{(\rho^{k}f, \rho^{l}g): (f, g)\in G_n^{\perp}, \,\|(f, g)\|_X=1\}.$$ 

\begin{lemma}\label{2410121927}
For $n\ge 1$, there exist two sequences $\rho_n>r_n>0$ such that 
\begin{itemize}
\item[\rm $(a)$] $a_n:=\underset{z\in  S_{r_n}^n}{\inf}\mathcal{J}(z)\to \infty$ as $n\to\infty$;
\item[\rm $(b)$] $b_n:=\underset{z\in Q^{2n+1}_{\rho_n}}{\sup}\mathcal{J}(z)\le 0$.
\end{itemize}
\end{lemma}
\begin{proof}
(a) For any $(\rho^{k}f, \rho^{l}g)\in S_{\rho}^n$ with $\rho\ge 1$, by Lemmas \ref{2211011024} and \ref{2408292319}, 
\[
\begin{aligned}
\mathcal{J}(\rho^{k}f, \rho^{l}g)
&\ge C_1 \rho^{k(1+1/p)}|f|_{1+\frac1p}^{1+\frac1p}+C_1\rho^{l(1+1/q)}|g|_{1+\frac1q}^{1+\frac1q}-\gamma_n\rho^{k+l}|f|_{1+\frac1p}|g|_{1+\frac1q}\\
& \ge C_2\rho^{\min\{k(1+1/p),\, l(1+1/q)\}}-\gamma_n\rho^{k+l}.
\end{aligned}
\]
If we let 
$$r_n:=\left(\frac{C_2}{2\gamma_n}\right)^{\frac{1}{k+l-\min\{k(1+1/p),\, l(1+1/q)\}}},$$
then for any $(r_n^{k}f, r_n^{l}g)\in S_{r_n}^n$, we have
\[
\mathcal{J}(r_n^{k}f, r_n^{l}g)\ge \left(\frac{C_2}{2}\right)^{\frac{k+l}{k+l-\min\{k(1+1/p),\, l(1+1/q)\}}}\gamma_n^{\frac{-\min\{k(1+1/p),\, l(1+1/q)\}}{k+l-\min\{k(1+1/p),\, l(1+1/q)\}}}\to\infty.
\]

%For any $(\rho^{k}f, \rho^{l}g)\in B_{r_n}^n$,
%\[
%\begin{aligned}
%\mathcal{J}(\rho^{k}f, \rho^{l}g)
%&\ge C_1 \rho^{k(1+1/p)}|f|_{1+\frac1p}^{1+\frac1p}+C_1\rho^{l(1+1/q)}|g|_{1+\frac1q}^{1+\frac1q}-C_3\rho^{k+l}|f|_{1+\frac1p}|g|_{1+\frac1q}\\
%& \ge C_2\rho^{\max\{k(1+1/p),\, l(1+1/q)\}}-C_3\rho^{k+l}.
%\end{aligned}
%\]
%Since $\rho\in [0, r_n]$ and $r_n\to 0$, we get $b_n=0$ for sufficiently large $n$.

\medskip
(b) For any $(\rho^k f, \rho^l g)\in Q^{2n+1}_{\rho}$, by Lemma \ref{2408300019}, we get as $\rho\to \infty$,
\[
\begin{aligned}
\mathcal{J}(\rho^{k}f, \rho^{l}g)
&\le C_3 \rho^{k(1+1/p)}|f|_{1+\frac1p}^{1+\frac1p}+C_3\rho^{l(1+1/q)}|g|_{1+\frac1q}^{1+\frac1q}-\rho^{k+l}\int_{\Omega}f\mathcal{A}g dx\\
& \le C_3 \rho^{k(1+1/p)}|f|_{1+\frac1p}^{1+\frac1p}+C_3\rho^{l(1+1/q)}|g|_{1+\frac1q}^{1+\frac1q}-C_4\rho^{k+l}\\
&\to -\infty.
\end{aligned}
\]
Since $F_{2n+1}$ has finite dimensions, the above limit  uniformly holds in $Q^{2n+1}_{\rho}$. Hence, we can find a sequence $\rho_n>r_n$ such that $(b)$ holds. The proof is done.
\end{proof}
\smallskip

Now we fix $\rho_n$ and $r_n$ satisfying Lemma \ref{2410121927}. Define
\[
\Gamma_n:=\{\gamma\in C(B_{\rho_n}^{2n+1}, X): \forall\; z\in B_{\rho_n}^{2n+1},\; \gamma(-z)=-\gamma(z);\;\gamma(z)=z \mbox{ for all } z\in Q_{\rho_n}^{2n+1}\}.
\]
Next, we give  the intersection property.
\begin{proposition}\label{2410121928}
For any $n\ge 1$, we have $\gamma(B_{\rho_n}^{2n+1})\cap S_{r_n}^n\neq \varnothing$ for all $\gamma\in \Gamma_n$.
\end{proposition}
\begin{proof}
We assume that  there exists some $\gamma\in \Gamma_n$ such that $\gamma(B_{\rho_n}^{2n+1})\cap S_{r_n}^n= \varnothing$.
We define the map $\mathcal{B}: [0, \rho_n]\times Q_1^{2n+1}\to X$ as 
\[
\mathcal{B}(t, (f, g)):=\gamma(t^{k}f, t^lg), \quad \forall \; t\in [0, \rho_n],\; (f, g)\in Q_1^{2n+1}.
\]
Let us denote
\[
X_1:=\{(\rho^{k}f, \rho^{l}g): (f, g)\in X, \,\|(f, g)\|_X=1,\; 0\le \rho\le r_n\}
\]
and
\[
X_2:=\{(\rho^{k}f, \rho^{l}g): (f, g)\in X, \,\|(f, g)\|_X=1, \; r_n\le \rho\}.
\]
Thus, we have 
\[
\mathcal{B}(0, Q_{1}^{2n+1})\subset X_1 \quad \mbox{and}\quad  \mathcal{B}(\rho_n, Q_{1}^{2n+1})\subset X_2.
\]

Firstly, we claim $X=X_1\cup X_2$. Indeed, for any $(f, g)\in X$, if $(f, g)=0$, it follows that $(f, g)\in X_1$. Therefore, we can simply assume $(f, g)\neq 0$. Let us set
\[
\widetilde{f}_t:=t^{-k}f, \quad \widetilde{g}_t:=t^{-l}g.
\]
Noting that 
\[
\underset{t\to \infty}{\lim}\|(\widetilde{f}_t, \widetilde{g}_t)\|_{X}=0,\quad \underset{t\to 0^+}{\lim}\|(\widetilde{f}_t, \widetilde{g}_t)\|_{X}=\infty,
\]
we then get that there exists some $t_0>0$ such that $\|(\widetilde{f}_{t_0}, \widetilde{g}_{t_0})\|_{X}=1$, hence 
$$(f, g)=(t_{0}^k\widetilde{f}_{t_0}, t_0^l\widetilde{g}_{t_0})\in X_1\cup X_2.$$
So $X\subset X_1\cup X_2$, and the claim holds.
\smallskip

Now we shall prove
\begin{equation}\label{2410222010}
{\rm Ind}({\rm Im}(\mathcal{B})\cap X_1\cap X_2)\le 2n.
\end{equation}
In fact, since $\gamma(B_{\rho_n}^{2n+1})\cap S_{r_n}^n= \varnothing$, we conclude that $A_n:={\rm Im}(\mathcal{B})\cap X_1\cap X_2\subset X\backslash \widetilde{A}_n$, where
\[
\widetilde{A}_n:=\{(\rho^k f, \rho^l g): (f, g)\in G_n^{\perp},\; \|(f, g)\|_X=1,\;\rho\ge 0\}.
\]
 Thus, for any $(r_n^kf, r_n^lg)\in A_n$ with $\|(f, g)\|_X=1$, we have $(\widetilde{P}_n f, P_n g)\neq 0$. Putting 
 $$B_n:=\{(\widetilde{P}_n f, P_n g): (r_n^kf, r_n^lg)\in A_n, \, \|(f, g)\|_X=1\},$$
 according to Lemma \ref{2410222003}(ii)(iii)(v), $B_n\subset G_n$ and Lemma \ref{2410140905}, we have
 \[
 {\rm Ind}(A_n)\le {\rm Ind}(B_n)\le {\rm dim}(G_n)=2n.
 \]
So \eqref{2410222010} holds. According to Lemma \ref{2410222003}(iv)(v) and \eqref{2410222010}, we have
\[
2n\ge {\rm Ind}({\rm Im}(\mathcal{B})\cap X_1\cap X_2)\ge {\rm Ind}(Q_1^{2n+1})=2n+1.
\]
That is a contradiction. The proof is done.
\end{proof}
\subsection{Proof of Theorem \ref{thm1}}
Under the assumptions $(H1)$ and $(H2)$, it follows from Lemma \ref{22112307} that $\mathcal{J}$ is of class $C^1$. Lemma \ref{2410222032} and $(H4)$ can deduce that $\mathcal{J}$ is even. In view of Proposition \ref{2410121928}, we can define the minimax levels
\[
M_{n}:=\underset{\gamma\in \Gamma_n}{\inf}\underset{z\in B_{\rho_n}^{2n+1}}{\max}\mathcal{J}(\gamma(z)).
\]
By Lemma \ref{2410121927} and Proposition \ref{2410121928}, we obtain $M_n\to\infty$ as $n\to\infty$. Using $(H3)$,  the Palais-Smale condition holds, which can refer to \cite[Lemma 4.8]{Z_arXiv2024}, so $M_n$ is a critical level. Thus, we derive infinitely many critical points. Finally, using Lemmas \ref{2410221530}, \ref{2401052112}, the proof is finished.

\section{Sublinear case: infinitely many solutions with negative energies}\label{2411022316}
In this section, we assume that $p,\, q$ are sublinear, i.e. $p,\, q>0$ and \eqref{2410222037} holds. For this situation, we can find $k,\, l>1$ such that
\begin{equation}\label{2211012223}
\frac{p}{p+1}<\frac{k}{k+l}, \quad \frac{q}{q+1}<\frac{l}{k+l}.
\end{equation}
For $m\ge n$, we  denote 
\[
B^{n, m}_\rho:=\{(t^{k}f, t^{l}g): (f, g)\in G_n^m, \,\|(f, g)\|_X= 1,\; 0\le t\le \rho\}
\]
and
\[
D^{n, m}_\rho:=\{(\rho^{k}f, \rho^{l}g): (f, g)\in G_n^m, \,\|(f, g)\|_X= 1\}.
\]

\begin{lemma}\label{2410131504}
There exists some $n_0\ge 1$ such that for any $n\ge n_0$, one can find two sequences $\rho_n>r_n>0$ such that $\rho_n\to 0$ as $n\to\infty$ and
\begin{itemize}
\item[\rm $(a)$] $\widetilde{a}_n:=\underset{z\in  S_{\rho_n}^{n-1}}{\inf}\mathcal{J}(z)>0$;
\item[\rm $(b)$] $\widetilde{b}_n:=\underset{z\in Q^{2n}_{r_n}}{\sup}\mathcal{J}(z)< 0$;
\item[\rm $(c)$] $\widetilde{d}_n:=\underset{z\in B_{\rho_n}}{\inf}\mathcal{J}(z)\to 0$ as $n\to\infty$, where  $B_{\rho}:=\{(t^{k}f, t^{l}g): (f, g)\in X, \,\|(f, g)\|_X= 1,\; 0\le t\le \rho\}$.
\end{itemize}
\end{lemma}
\begin{proof}
(a) For any $(\rho^{k}f, \rho^{l}g)\in S_{\rho}^{n-1}$ with $0<\rho\le 1$, by Lemmas \ref{2211011024}, \ref{2408292319},
\[
\begin{aligned}
\mathcal{J}(\rho^{k}f, \rho^{l}g)
&\ge C_1 \rho^{k(1+1/p)}|f|_{1+\frac1p}^{1+\frac1p}+C_1\rho^{l(1+1/q)}|g|_{1+\frac1q}^{1+\frac1q}-\gamma_{n-1}\rho^{k+l}|f|_{1+\frac1p}|g|_{1+\frac1q}\\
& \ge C_2\rho^{\max\{k(1+1/p),\, l(1+1/q)\}}-\gamma_{n-1}\rho^{k+l}.
\end{aligned}
\]
 When $n$ is large enough, we can define
$$\rho_n:=\left(\frac{2\gamma_{n-1}}{C_2}\right)^{\frac{1}{\max\{k(1+1/p),\, l(1+1/q)\}-(k+l)}}<1,$$
then for any $(\rho_n^{k}f, \rho_n^{l}g)\in S_{\rho_n}^{n-1}$, we have
\[
\mathcal{J}(\rho_n^{k}f, \rho_n^{l}g)\ge \left(\frac{2}{C_2}\right)^{\frac{k+l}{\max\{k(1+1/p),\, l(1+1/q)\}-(k+l)}}\gamma_{n-1}^{\frac{\max\{k(1+1/p),\, l(1+1/q)\}}{\max\{k(1+1/p),\, l(1+1/q)\}-(k+l)}}>0.
\]

(b) For any $(f, g)\in F_{2n}$ with $\|(f, g)\|_X=1$, by Lemma \ref{2408300019}, we get that for $\rho\in (0, 1)$,
\[
\begin{aligned}
\mathcal{J}(\rho^{k}f, \rho^{l}g)
&\le C_3 \rho^{k(1+1/p)}|f|_{1+\frac1p}^{1+\frac1p}+C_3\rho^{l(1+1/q)}|g|_{1+\frac1q}^{1+\frac1q}-\rho^{k+l}\int_{\Omega}f\mathcal{A}g dx\\
& \le C_4 \rho^{\min\{k(1+1/p), l(1+1/q)\}}-C_5\rho^{k+l}.
\end{aligned}
\]
By $0<k+l<\min\{k(1+1/p), l(1+1/q)\}$, there exists $r_n\in (0, \rho_n)$ such that $\mathcal{J}(r_n^{k}f, r_n^{l}g)<0$.

\medskip
(c) For any $(f,g)\in X$ with $\|(f, g)\|_X= 1$, it holds
$$
\Big|\int_{\Omega}g\mathcal{A}f dx \Big| \leq |g|_{1+\frac1q}|\mathcal{A}f|_{1+q} \leq C|g|_{1+\frac1q}|\mathcal{A}f|_{W^{2, 1+\frac1p}(\Omega)} \leq C|g|_{1+\frac1q}|f|_{1+\frac1p} \leq C.
$$
So,  for any $(t^{k}f, t^{l}g)\in B_{\rho_n}$, $0\le t\le \rho_n$,
$$
-Ct^{k+l} \leq \mathcal{J}(t^{k}f, t^{l}g) \leq C \Big(  t^{k(1+1/p)}|f|_{1+\frac1p}^{1+\frac1p}+t^{l(1+1/q)}|g|_{1+\frac1q}^{1+\frac1q}+t^{k+l}\Big). 
$$
Since $\rho_n\to 0$ as $n\to\infty$, we conclude that
$\widetilde{d}_n\to 0$ as $n\to\infty$.
\end{proof}

Next, let $r_n$ and $\rho_n$ be given in Lemma \ref{2410131504}. For $m\ge 2n$, let us denote
\[
\Gamma_n^m:=\{\gamma\in C(B_{\rho_n}^{n, m}, G_m):  \gamma(-z)=-\gamma(z),\, \forall\, z\in B_{\rho_n}^{n, m};\;\gamma(z)=z,\, \forall\, z\in D_{\rho_n}^{n, m}\}.
\]
Similar to Proposition \ref{2410121928}, we have the following intersection property.
\begin{proposition}\label{2410132041}
For $m\ge 2n$, $\gamma(B_{\rho_n}^{n, m})\cap Q_{r_n}^{2n}\neq \varnothing$ for all $\gamma\in \Gamma_n^m$.
\end{proposition}
\begin{proof}
If there is a $\gamma\in \Gamma_n^m$ satisfying $\gamma(B_{\rho_n}^{n, m})\cap Q_{r_n}^{2n}= \varnothing$,
we can define the map $\mathcal{B}: [0, \rho_n]\times D_{1}^{n, m}\to G_m$ as 
\[
\mathcal{B}(t, (f, g)):=\gamma(t^{k}f, t^lg), \quad \forall \; t\in [0, \rho_n],\; (f, g)\in D_{1}^{n, m}.
\]
We denote
\[
Y_1:=\{(\rho^{k}f, \rho^{l}g): (f, g)\in G_m, \,\|(f, g)\|_X=1,\; 0\le \rho\le r_n\}
\]
and
\[
Y_2:=\{(\rho^{k}f, \rho^{l}g): (f, g)\in G_m, \,\|(f, g)\|_X=1, \; r_n\le \rho\}.
\]
As we proved in Proposition \ref{2410121928}, there hold that
\[\mbox{
$G_m=Y_1\cup Y_2$,\quad $\mathcal{B}(0, D_{1}^{n, m})\subset Y_1$\quad  and\quad  $\mathcal{B}(\rho_n, D_{1}^{n, m})\subset Y_2$.}
\]

Now we claim that 
\begin{equation}\label{2410222006}
{\rm Ind}(A_{m, n})\le 2m-2n,
\end{equation}
where $A_{m,n}:={\rm Im}(\mathcal{B})\cap Y_1\cap Y_2$. Since  $\gamma(B_{\rho_n}^{n, m})\cap Q_{r_n}^{2n}= \varnothing$, we conclude that $A_{m,n}\subset G_m\backslash \widetilde{A}_{m,n}$, where
\[
\widetilde{A}_{m,n}:=\{(\rho^k f, \rho^l g): (f, g)\in F_{2n},\;\|(f, g)\|_X=1,\; \rho\ge 0\}.
\]

On the other hand, since $G_m$ has finite dimensions and $F_{2n}$ is a linear subspace of $G_m$, there are a subspace $L_{2m-2n}$ of $G_m$ and a linear projection operator $\overline{P}: G_m\to F_{2n}$ such that
\[
G_m=F_{2n}\oplus L_{2m-2n},
\]
and for any $(f, g)\in G_m$, there exists a unique decomposition
\[
\overline{P}(f, g)\in F_{2n},\quad (f, g)-\overline{P}(f, g)\in L_{2m-2n}.
\]

Finally, for any $(r_n^kf, r_n^lg)\in A_{m,n}$ with $\|(f, g)\|_X=1$, we put
 $$B_{m,n}=\{(f,g)-\overline{P}(f, g): (r_n^kf, r_n^lg)\in A_{m,n}, \, \|(f, g)\|_X=1\}.$$
According to Lemma \ref{2410222003}(ii)(iii)(iv) and $B_{m,n}\subset L_{2m-2n}$,
 \[
 {\rm Ind}(A_{m,n})\le {\rm Ind}(B_{m,n})\le {\rm dim}(L_{2m-2n})=2m-2n.
 \]
So the claim holds. According to Lemma \ref{2410222003}(iv)(v) and \eqref{2410222006}, we have
\[
2m-2n\ge {\rm Ind}({\rm Im}(\mathcal{B})\cap Y_1\cap Y_2)\ge {\rm Ind}(D_1^{n, m})=2m-2n+2.
\]
That is a contradiction. The proof is done.
\end{proof}

\subsection{Proof of Theorem \ref{thm2}}

Now for $m\ge 2n$, it follows from Lemma \ref{2410131504} and Proposition \ref{2410132041} that the following minimax levels are well-defined
\[
c_n^m:=\underset{\gamma\in \Gamma_n^m}{\sup}\underset{z\in B_{\rho_n}^{n, m}}{\min}\mathcal{J}(\gamma(z)).
\]
Lemma \ref{2410131504} and Proposition \ref{2410132041} can deduce that $\widetilde{d}_n\le c_n^m\le\widetilde{b}_n<0$ for all $m\ge 2n$. On the one hand, since $c_n^m$ is bounded, up to a subsequence, there exists $c_n<0$ such that
\[
c_n^m\to c_n,\quad\mbox{as }m\to\infty.
\]
Assume that $(w_{m, n}, y_{m, n})\in G_m$ is a critical point of $\mathcal{J}|_{G_m}$ corresponding to the level $c_n^m$. Since $\mathcal J$ is coercive, $\{(w_{m, n}, y_{m, n})\}$ is bounded in $X$. We can assume $(w_n, y_n)$ is a weak convergence point of $\{(w_{m, n}, y_{m, n})\}$ in $X$. Since $\mathcal{A}: L^{1+\frac{1}{p}}(\Omega)\to L^{q+1}(\Omega)$ is compact and $X=\cup_{m\ge 1}G_m$, it is easy to conclude $(w_n, y_n)$ is a critical point of $\mathcal{J}$. Thus $c_n$ is a critical value. By $\widetilde{d}_n\to 0$, we have $c_n\to 0^-$. As a result, we obtain infinitely many critical points of $\mathcal{J}$, whose energies are negative and tend to $0$. Finally, using Lemmas \ref{2410221530}, \ref{2401052112}, we finish the proof.

\bigskip

\noindent \textbf{Statements and Declarations} 
\medskip

\noindent \textbf{Conflict of interest} The authors declare that they have no conflict of interest.\\
\noindent \textbf{Data availability} Data sharing is not applicable to this article as no new data
were created or analyzed in this study.


\newpage
\begin{thebibliography}{99}
\bibitem{B1993} T. Bartsch, Infinitely many solutions of a symmetric Dirichlet problem, {\it Nonlinear Anal.} {\bf 20} (1993), 1205-1216. 

\bibitem{Bd1999} T. Bartsch and D.G. de Figueiredo, Infinitely many solutions of nonlinear elliptic systems, Topics in nonlinear analysis, {\it Progr. Nonlinear Differential Equations Appl.}, 35, Birkh\"auser, Basel, (1999), 51-67.

\bibitem{BW1995} T. Bartsch and M. Willem, On an elliptic equation with concave and convex nonlinearities, {\it Proc. Amer. Math. Soc.} {\bf 123} (1995),  3555-3561.

\bibitem{BR1979} V. Benci and P.H. Rabinowitz, Critical point theorems for indefinite functionals, {\it Invent. Math.} {\bf 52} (1979), 241-273.



\bibitem{CV1995} P.H. Cl\'ement and  R.C.A.M. van der Vorst, On a semilinear elliptic system, {\it Differential Integral Equations} {\bf 8} (1995), 13171-329.

\bibitem{ddR_JFA2005} D.G. de Figueiredo, J.M. do \'{O} and B. Ruf, An {O}rlicz-space approach to superlinear elliptic systems, {\it J. Funct. Anal.} {\bf 224} (2005), 471-496.

\bibitem{dF1994} D.G. de Figueiredo and P.L. Felmer, On superquadratic elliptic systems, {\it Trans. Amer. Math. Soc.} {\bf 343} (1994), 99-116.

%\bibitem{Ding_TMNA1997} Y. Ding, Infinitely many entire solutions of an elliptic system with symmetry, {\it Topol. Methods Nonlinear Anal.} {\bf 9} (1997), 313-323.


\bibitem{FR1978} E.R. Fadell and P.H. Rabinowitz, Generalized cohomological index theories for Lie
group actions with an application to bifurcation questions for Hamiltonian systems,
{\it Invent. Math.} {\bf 45} (1978), 139-174.

\bibitem{FZ2003} X.L. Fan and Q.H. Zhang, Existence of solutions for $p(x)$-Laplacian Dirichlet problem,
{\it Nonlinear Anal.} {\bf 52} (2003), 1843-1852. 


\bibitem{HMv1998} J. Hulshof, E. Mitidieri and R.C.A.M. van der Vorst, Strongly indefinite systems with critical Sobolev exponents, {\it Trans. Amer. Math. Soc.} {\bf 350} (1998), 2349-2365.


\bibitem{Hv1993} J. Hulshof and R.C.A.M. van der Vorst, Differential systems with strongly indefinite variational structure, {\it J. Funct. Anal.} {\bf 114} (1993), 32-58.

\bibitem{Liu2010} S. Liu, On superlinear problems without the Ambrosetti and Rabinowitz condition, {\it Nonlinear Anal.} {\bf 73} (2010), 788-795.


\bibitem{Mitidieri1993} E. Mitidieri, A Rellich type identity and applications: Identity and applications, {\it Comm. Partial Differential Equations} {\bf 18} (1993), 125-151.

\bibitem{PSY2016} K. Perera,  M. Squassina and  Y. Yang, Bifurcation and multiplicity results for critical fractional $p$-Laplacian problems,
{\it Math. Nachr.} {\bf 289} (2016), 332-342.

\bibitem{Struwe2008} M. Struwe, Variational methods, Applications to nonlinear partial differential equations and Hamiltonian systems, Fourth edition, Springer-Verlag, Berlin, 2008.


\bibitem{YZ_ANS2024} D. Ye  and W. Zhang,  Existence and multiplicity of solutions for fractional $p$-Laplacian equation involving critical concave-convex nonlinearities,
{\it Adv. Nonlinear Stud.} {\bf 24} (2024), 895-921.

\bibitem{Z_arXiv2024} W. Zhang, Variational method for fractional Hamiltonian system in bounded domain, arXiv:2404.00687 (2024). https://arxiv.org/abs/2404.00687



\end{thebibliography}
\end{document}